\documentclass[article]{IEEEtran}

\usepackage{amsmath, amssymb, amsfonts} 
\usepackage{algorithm}
\usepackage{algorithmic}
\usepackage{array}
\usepackage{subfig}

\usepackage{textcomp}
\usepackage{stfloats} 
\usepackage{url}
\usepackage{verbatim}
\usepackage{graphicx}
\usepackage{cite} 
\usepackage{listings} 
\usepackage{booktabs} 
\usepackage{comment}
\usepackage[table,xcdraw]{xcolor} 
\usepackage{adjustbox} 
\usepackage{wrapfig}

\begin{document}

\title{Attention is All You Need to Optimize \\ Wind Farm Operations and Maintenance}

\author{Iman Kazemian$^\dagger$, Murat Yildirim$^\dagger$ and Paritosh Ramanan$^\ddagger$\\
			$^\dagger$ Industrial and Systems Engineering, Wayne State University, Detroit, US \\
             $^\ddagger$ Industrial Engineering and Management, Oklahoma State University, Stillwater, US }



\maketitle

\begin{abstract}
Operations and maintenance (O\&M) is a fundamental problem in wind energy systems with far reaching implications for reliability and profitability. Optimizing O\&M is a multi-faceted decision optimization problem that requires a careful balancing act across turbine level failure risks, operational revenues, and maintenance crew logistics. The resulting O\&M problems are typically solved using large-scale mixed integer programming (MIP) models, which yield computationally challenging problems that require either long-solution times, or heuristics to reach a solution. To address this problem, we introduce a novel decision-making framework for wind farm O\&M that builds on a multi-head attention (MHA) models, an emerging artificial intelligence methods that are specifically designed to learn in rich and complex problem settings. The development of proposed MHA framework incorporates a number of modeling innovations that allows explicit embedding of MIP models within an MHA structure. The proposed MHA model (i) significantly reduces the solution time from hours to seconds, (ii) guarantees feasibility of the proposed solutions considering complex constraints that are omnipresent in wind farm O\&M, (iii) results in significant solution quality compared to the conventional MIP formulations, and (iv) exhibits significant transfer learning capability across different problem settings. 
\end{abstract}


\section{Introduction}
Operations and maintenance (O\&M) is a central challenge and a catalyst for increasing competitiveness of offshore wind energy. As the U.S. targets the installation of 30 GW of offshore wind energy by 2030 \cite{office2022fact}, the cost of maintaining these offshore installations is becoming increasingly significant. Globally, the offshore wind industry, with over 59 GW of installed capacity across 293 operating projects, incurs substantial O\&M costs, estimated at around \$3 billion annually\cite{office2022fact}. These costs contribute significantly to the levelized cost of energy (LCOE) for offshore wind, with O\&M expenses alone making up around 25\% of LCOE in offshore wind farms \cite{us2022levelized}. Evidently, an effective O\&M framework can have overarching implications for increasing operational availability, effectiveness and reducing the costs associated with maintenance and related logistics; hence becoming a catalyst for ensuring long-term competitiveness of offshore wind energy systems.


O\&M is a multifaceted problem incorporating complex interdependencies across operations and maintenance. Conventional models focus primarily on the operational aspect of O\&M. These models optimize operational goals with the assumption that a periodic maintenance (e.g. annual maintenance) would ensure an acceptable level of turbine availability. In reality, wind turbine failures are frequent, and accurately predicting these failures offers substantial potential for enhancing both operational efficiency and maintenance strategies. Recently, there has been a significant focus on a set of models that harness real-time sensor data from wind turbines to better predict wind turbine conditions and remaining life distributions, a field called prognostics \cite{rezamand2020critical, alvarez2022power}. The assumption is that improved prognostic predictions on turbine lifetime can lead to better O\&M outcomes. However, for complex problems like wind farm O\&M, prognostic predictions constitute only half the battle. A significant challenge remains in translating these predictions into effective strategies for managing the complex and interconnected O\&M decisions that accompany them.

Research on prognostics-driven O\&M has three research streams. The first stream focuses on using sensor data to optimize maintenance timing for individual components or turbines, aiming to balance premature maintenance—which can waste equipment life—and delayed maintenance, which raises the risk of failure. Examples include incorporating failure probabilities into maintenance decisions \cite{raza2019optimal}, accounting for dynamic weather conditions \cite{byon2013wind}, and optimally replacing degrading components in partially observable environments \cite{flory2015optimal}. The second stream models multiple components or turbines by modeling on the economic benefits of grouping maintenance actions, known as opportunistic maintenance. This includes condition-based maintenance models using multi-threshold O\&M\cite{tian2011condition}, combining CBM with opportunistic thresholds\cite{zhou2019opportunistic}, and managing spare parts supply in wind farms \cite{zhu2022joint}. The third research stream models farm-level interactions, extending beyond opportunistic maintenance. These models require complex optimization frameworks that integrate prognostic predictions into large-scale decision-making processes, such as capturing dependencies across turbines\cite{yildirim2017integrated}, applying stochastic optimization with chance constraints\cite{fallahi2022chance}, considering accessibility impacts\cite{papadopoulos2021seizing}, and scheduling maintenance for multi-component turbines\cite{fallahi2022chance}. While these models provide robust frameworks for optimizing O\&M at various levels, they are often computationally intensive and difficult to solve as wind farm scale and complexities increase.

To tackle the computational challenges, significant attention has been directed toward leveraging machine learning (ML) models to either accelerate or replace traditional decision optimization models. Acceleration has been achieved by developing ML models that approximate lower bounds \cite{baltean2018selecting}, or enhance branching policies \cite{applegate1998solution,alvarez2017machine,gasse2019exact,he2014learning,lodi2017learning} during the optimization process. Concurrently, another research area has focused on creating ML algorithms that directly output optimal decisions, particularly excelling in well-structured problems with minimal constraints, such as the Traveling Salesman problem and its variants, using techniques like pointer networks\cite{vinyals2015pointer}, deep reinforcement learning\cite{kool2018attention}, and graph neural networks \cite{khalil2017learning,khalil2016learning}. However, for more general class of constrained optimization problems, such as O\&M scheduling, current ML-based models often face difficulties with feasibility, scalability, and generalization across different instances.

Attention transformation methods have recently been investigated to propose ML based methods that can inherently capture complex inter-dependencies in complex optimization models. Since the introduction of the attention mechanism in the seminal work ``Attention is All You Need" by Vaswani et al. \cite{vaswani2017attention}, this approach demonstrated significant promise in a wide range of applications, including Unmanned Aerial Vehicle (UAV) scheduling in humanitarian response \cite{wan2024deep}, real-time human-robot collaboration \cite{wang2019learning}, maximal covering location problem\cite{wang2023deepmclp}, crane scheduling in storage yards \cite{jin2023deep}, and task management in manufacturing \cite{wang2022solving}. Within the attention transformation methods, multi-head attention transformer-based approaches are well-suited for addressing permutation-based problems such as Travelling Salesman Problems (TSPs) and Vehicle Routing Problems (VRPs) \cite{kool2018attention} due to their sequential encoding and decoding frameworks, and were shown to achieve superior performance compared to other approaches \cite{liu2022good}. These approaches can effectively capture the order of node selection, which is crucial in optimizing routing.

To date, researchers have successfully employed attention mechanisms to solve well-defined optimization problems characterized by relatively simple constraints and reward functions. For example, in the context of the Traveling Salesman Problem (TSP), modeled as a sequential selection problem, the task reduces to choosing the next visit location from a subset of nodes \cite{kool2018attention}. However, there hasn't been an application of multi-head attention networks to more complex problems, such as those encountered in the decision making within energy sector. These problems involve challenging features, including time-varying maintenance costs, unexpected failures, fluctuating market prices, crew visit costs, logistical considerations, and opportunistic maintenance interactions. Transformer-based approaches hold significant potential to assist wind farm operators in identifying optimal solutions efficiently, thereby enhancing maintenance scheduling and resource allocation.

In this work, we propose a multi-head attention framework model to optimize the wind farm O\&M scheduling problem, which we refer to as \textit{AttenCOpt: Attention-based constrained optimization model}. Our approach integrates sensor-driven prognostic predictions from a large number of turbines, with optimal O\&M decision making considering factors such as opportunistic maintenance, failure scenarios, degradation signals, and the number of maintenance activities per period. We incorporate production management decisions that account for market prices and maximum production limit scenarios over time. We propose a novel method that guarantees feasibility through use of model structure and mapping of constraints to the masking procedures. Transformers, with their self-attention mechanisms, excel in handling large and complex datasets by dynamically weighting the importance of different input elements. This capability makes attention transformation a powerful approach for solving wind farm maintenance and operation problems, as it enhances the model's ability to understand and leverage the structure of the optimization problem, leading to improved solution quality and computational efficiency. 
The proposed approach enhances the robustness of the solution process and increases the practical applicability of ML models in solving complex OR problems. Our contributions are as follows:

\begin{enumerate}
    \item We develop a novel multi-head attention network (MHA) model that integrates (i) sensor-driven prognostic predictions for wind turbine remaining life with (ii) optimal O\&M scheduling decisions for large-scale wind farms. The proposed model architecture enhances the model’s capability to process and learn from complex, high-dimensional data, improving the quality of scheduling decisions. 

    \item Feasibility and optimality of the model is attained by developing generalizable mappings from constraints and objective functions of the mathematical programming formulations into the multi-head attention networks architecture. Innovations in the proposed model include:
    \begin{itemize}
    \item (i) Masking procedures that are used to build a one-to-one mapping across optimization constraints to their projections in conditional masking procedures. 
    \item (ii) Novel reward functions to model the objective function by combining production revenues, and maintenance \& logistical expenses, to offer a comprehensive evaluation of O\&M decisions. 
    \end{itemize}

    \item We offer a cohesive, end-to-end framework that interweaves multi-head attention framework and optimization techniques during training and testing to produce efficient and effective maintenance schedules.

    \item  We demonstrate the effectiveness of the proposed models in knowledge transfer across varying problem settings by conducting training and testing on different scenarios. Notably, the models achieve significant success, when the test case study is smaller in scale (i.e. number of wind turbines) compared to the training case study.

\end{enumerate}

We demonstrate the effectiveness of the proposed model through comprehensive set of experiments conducted using real-world degradation datasets. The remainder of the paper proceeds as follows: Section 2 provides a detailed Problem Description, outlining the specific characteristics mathematical model of wind farm maintenance scheduling problem. Section 3 presents the transformation of the Wind Farm Maintenance Scheduling Problem as a Learning Problem, where we define the key inputs, outputs, and features necessary for training the model. In Section 4, we introduce the Attention Model, describing its architecture, the role of multi-head attention, and how it captures temporal and non-temporal features in the context of turbine maintenance. Finally, Section 5 discusses the Experiments conducted to evaluate the proposed model.

\section{Problem Description}

\begin{figure*}[htbp]
	\centering
	\includegraphics[width=1\linewidth]{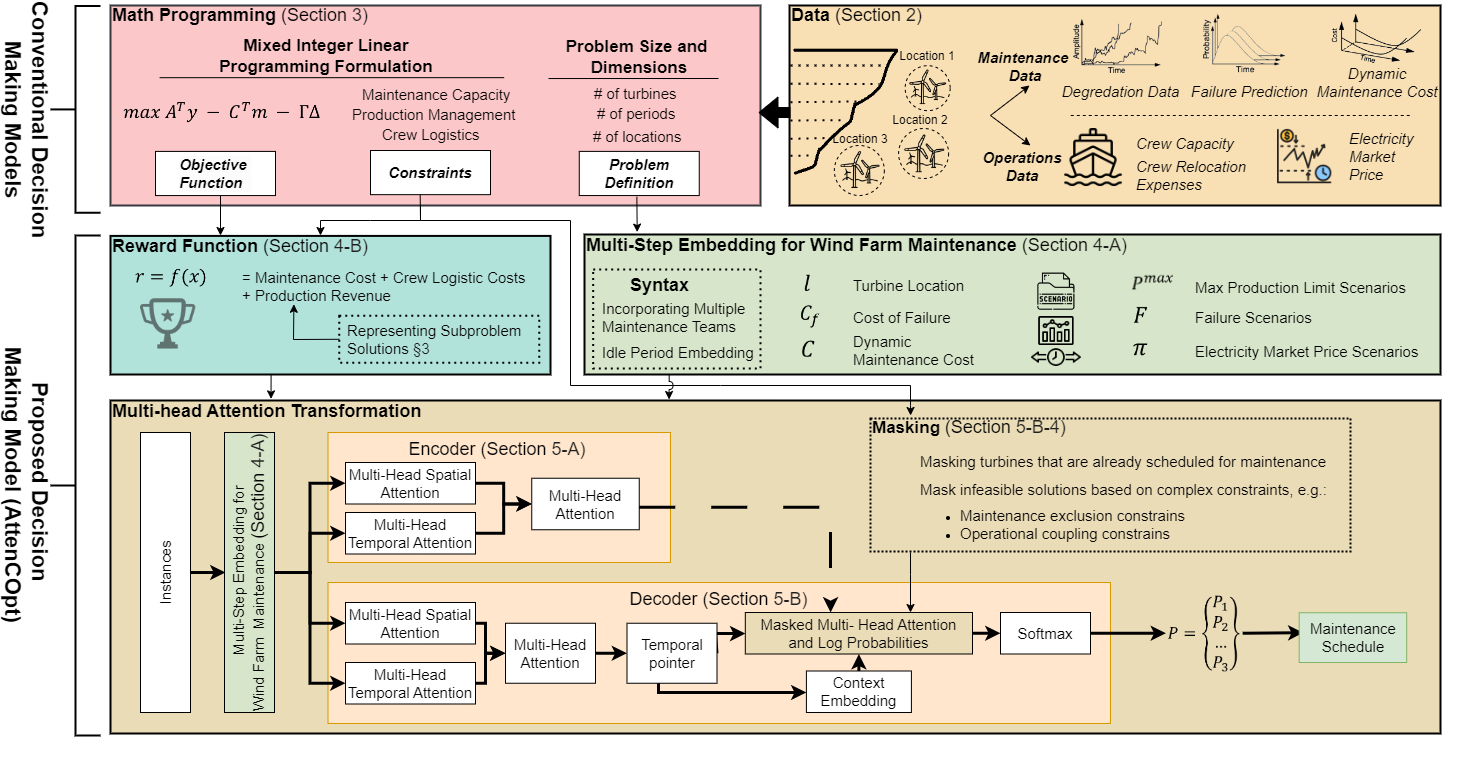} 
	\caption{Transitioning from data and math programming to AttenCOpt}
	\label{fig: Main_fig}
\end{figure*}

We study the fundamental problem of O\&M scheduling in wind energy systems, where the focus is to maximize operational profit from a fleet of wind turbines subject to numerous constrains related to operations (e.g. production of wind turbines as a function of availability and wind power) and maintenance (e.g. routing of maintenance crew). The proposed model is cast in a two-stage optimization formulation, where the first stage variables include the routing of maintenance team, and maintenance decisions. Customary to two-stage optimization models, we assume that these first set of decisions are made prior to the realization of uncertainty. Following the uncertainty realization, the production variables are assigned to model the resulting production from the wind turbines.

The proposed model considers multiple sources of uncertainty. The first two sets of uncertainty involves operational factors like market prices and wind power realization during each period. The second type pertains to degradation, which governs how the turbine degrades and fails over time. Both types of uncertainties are integrated into the model through use of scenarios. In line with the sources of uncertainty, we also integrate two types of data sources within a decision making framework. The first type deals with operations data, such as predictions on future price and wind power, logistics costs and parameters associated with crew routing, among others. The second data type deals with maintenance, which includes sensor information from the turbines and the associated predictions on remaining life and maintenance costs. The proposed model uses sensor-driven prognostic predictions to gain improved visibility on lifetime of the turbines, and leverage these predictions to optimize maintenance that balances (i) early maintenances which reduce equipment lifetime, with (ii) late maintenances which increases failure risks. The sensor-driven prognostic framework and the associated maintenance cost predictions have been outlined in  \cite{yildirim2017integrated,yildirim2016sensor1,yildirim2016sensor2}. Some of this data changes dynamically across the considered uncertainty scenarios, and time periods, enabling the proposed model to model dynamically evolving conditions.

To reflect the realistic operational conditions in wind farm O\&M, we model a fleet of turbines that are geographically distributed across multiple locations. Each time the maintenance crew visits a new location, a fixed visit cost is incurred. 
Evidently, we model a central decision-tradeoff for wind farm O\&M that encourages the maintenance teams to group turbine maintenances together in an effort to reduce the incurred visit costs. This tradeoff, called \textit{opportunistic maintenance}, is particularly important for offshore wind farms, where each visit to a wind farm introduces significant logistics costs associated with vessel and crew travel. Our objective also combines maintenance costs, which are functions of degradation signals and scheduling decisions, with revenues from energy production, adjusted across all scenarios based on actual market prices and turbine operational statuses as well as the crew visit cost. This approach allows for flexible scheduling that can adapt to varying operational conditions and degradation signals, optimizing both maintenance expenditures and energy production revenue.

As highlighted in Figure \ref{fig: Main_fig}, we solve this fundamental problem through two methods. The first method is the conventional benchmark method which uses mathematical programming model. In the second method, we introduce our proposed decision-making model, AttenCOpt, which uses an MHA framework to solve the same problem with feasibility guarantees. To this end, we showcase how to explicitly embed the objective function, constraints and problem definition of the mathematical programming formulation, within the proposed AttenCOpt. Figure \ref{fig: Main_fig} provides a holistic overview of this transition, and showcases how the dependencies are embedded across multiple components of the proposed MHA framework. Sections that elucidate the associated model developments are also highlighted in the figure \ref{fig: Main_fig}. In the next sections, we will first introduce the model as a conventional mixed integer program (MIP), and then develop the associated MHA framework.

\section{Wind Farm Maintenance Scheduling Problem as a Mathematical Formulation}

As introduced in Section II, we model the proposed wind farm O\&M model using a MIP formulation. The associated decision variables and parameters are provided below.

\subsection{Decision Variables}
\begin{itemize}
     \item \( g_{j,t}\):= 1, if maintenance team is in location \(j\) at period \( t \); otherwise 0.
         \item \( \delta_t \):= 1, if the location of maintenance team changes at period \( t \); otherwise 0.
             \item \( m_{i,t} \):= 1, if maintenance scheduled for turbine \(i\) at period \( t \); otherwise 0
\end{itemize}
\begin{itemize}
    \item \( y^s_{i,t} \): Production level for turbine \( i \) in period \( t \) under scenario \( s \). This variable reflects the responsive action taken once the actual conditions (such as market prices and turbine operational status) are known.

\end{itemize}

\subsection{Parameters}

\begin{itemize}
    \item \( I \):= Total number of turbines.
    \item \( T \):= Total number of periods.
    \item \( M \):= Maintenance capacity in each period.
    \item \( C_{i,t} \):= Estimated maintenance cost for turbine \( i \) in period \( t \).
    \item \( C^f \):= Estimated cost of failure.
    \item \( \Delta \):= The cost of changing location for maintenance team
    \item \( L_{i,j} \):= 1, if turbine \( i \) is located at position \( j \); otherwise 0.
\end{itemize}
\begin{itemize}
    \item \( S \):= Number of scenarios representing various possible future states.
    \item \( \pi^s_{t} \):= Realized price per unit of production in period \( t \) under scenario \( s \).
    \item \( P^s_{i,t} \):= Maximum production rate of turbine \( i \) in period \( t \) under scenario \( s \).
    \item  \( F^s_{i} \):= Failure time of turbine \( i \) under scenario \( s \).

\end{itemize}

\subsection{Objective Function}
The objective function maximizes the expected net profit by considering multiple revenue and cost factors.  
The first term models the revenue generated from production, which is a function of electricity price and production. Production itself revolves as a function of wind power and turbine availability. The second and third terms capture the preventive and corrective maintenance costs, which are determined using the scenario-specific failure time realizations, $F_i^s$. Finally, the fourth term models the visit costs by penalizing for every time a maintenance team changes location. 

\begin{equation}
\label{eq:originalobj}
\begin{aligned}
\text{Maximize} \quad \frac{1}{|S|} \biggl[ &\sum_{s=1}^{S} \sum_{i=1}^{I} \sum_{t=1}^{T}  \left(\pi^s_{t} \cdot y^s_{i,t} \right) - \\
&\biggl[ \sum_{s=1}^{S} \sum_{i=1}^{I} \sum_{t=1}^{F^s_i-1}  \left(C_{it} \cdot m_{i,t}\right) \\
&+ \sum_{s=1}^{S} \sum_{i=1}^{I} \sum_{t=F^s_i}^{T} \left( C^f \cdot m_{i,t} \right) \biggr] \biggr] - \sum_{t=1}^T \delta_t \cdot \Delta
\end{aligned}
\end{equation}

\subsection{Constraints}
\subsubsection{Production Management Constraints}
Constraints \eqref{eq:production_constraint} and \eqref{eq:ailalibility_constraint2} ensure that the production level \( y^s_{i,t} \) for turbine \( i \) in period \( t \) under scenario \( s \) cannot exceed the maximum production rate \( R^s_{i,t} \) if that turbine is available in the same period. When unavailable, the turbine production is set to zero.

\begin{equation}
\label{eq:production_constraint}
y^s_{i,t} \leq P^s_{i,t} (1-m_{i,t})  \quad \forall i \in I, \, t \leq F_{i}^s-1, \, s \in S
\end{equation}

\begin{equation}
\label{eq:ailalibility_constraint2}
y^s_{i,t} \leq P^s_{i,t} \sum_{l=1}^{t-1} m_{i,l} \quad \forall i \in I, \, t \ge F_{i}^s, \, s \in S
\end{equation}

\subsubsection{Maintenance Capacity Constraints}

Equation~\ref{eq:Maintenance Scheduling Constraints} ensures that the total number of turbines scheduled for maintenance in any given period \( t \) does not exceed the maintenance capacity \( M \).

\begin{equation}
\label{eq:Maintenance Scheduling Constraints}
\sum_{i=1}^{I} m_{i,t} \leq M \quad \forall t \in T
\end{equation}

Equation~\ref{eq:Maintenance Location Consistency} enforces that the maintenance decision \( m_{i,t} \) for turbine \( i \) in period \( t \) at location \( j \) can be conducted only if the maintenance team is in location \( g_{j,t} \):

\begin{equation}
\label{eq:Maintenance Location Consistency}
m_{i,t} \cdot L_{i,j} \leq g_{j,t} \quad \forall i \in I, \, t \in T, \, j \in J
\end{equation}

Equation ~\ref{eq:Maintenance Constraints} ensures that each turbine \( i \) is scheduled for maintenance exactly once during the planning horizon \( T \). 

\begin{equation}
\label{eq:Maintenance Constraints}
\sum_{t=1}^{T} m_{i,t} = 1 \quad \forall i \in I
\end{equation}

\subsubsection{Crew Logistics Constraints}
Equation~\ref{eq:location changed} sets the binary variable \( \delta_t \) to 1 if the location of the maintenance team changes from period \( t-1 \) to period \( t \):

\begin{equation}
\label{eq:location changed}
 g_{j, t-1} \leq  g_{j, t} + \delta_t \quad \forall t \in \{2, \ldots, T\}, \, j \in J
\end{equation}

\subsubsection{Non-negativity and Integrity Constraints}
All decision variables must adhere to non-negativity and integrity constraints:
\begin{equation}
g_{j,t} \in \{0, 1\} \quad \forall j \in J, \, t \in T
\end{equation}

\begin{equation}
m_{i,t} \in \{0, 1\} \quad \forall i \in I, \, t \in T
\end{equation}
\begin{equation}
\label{eq:domain}
\delta_{t} \in \{0, 1\} \quad \forall t \in T
\end{equation}

\textbf{Lemma:} The objective function~\ref{eq:originalobj} can be reformulated as:

\begin{equation}
\begin{aligned}
\label{obj}
\text{Minimize} \quad \frac{1}{|S|} \biggl[ &\sum_{s=1}^{S} \sum_{i=1}^{I} \sum_{t=1}^{F^s_i-1} \left( C_{i,t} + \pi^s_{t} \cdot P^s_{i,t} \right) \cdot m_{i,t} + \\
&\sum_{s=1}^{S} \sum_{i=1}^{I} \sum_{t=F^s_i}^{T} \left( C^f + \sum_{l=F^s_i}^t \pi^s_{l} \cdot P^s_{i,l} \right) \cdot m_{i,t} \biggr] \\
&+ \sum_{t=1}^T \delta_t \cdot \Delta.\\
&\text{S.t:}~(
\ref{eq:Maintenance Scheduling Constraints}-\ref{eq:domain})
\end{aligned}
\end{equation}
The proof is provided in Appendix \ref{appendix:proof1}.

\section{Wind Farm Maintenance Scheduling Problem as a Learning Problem}
\begin{figure*}[h]
    \centering
    \includegraphics[width=1\linewidth]{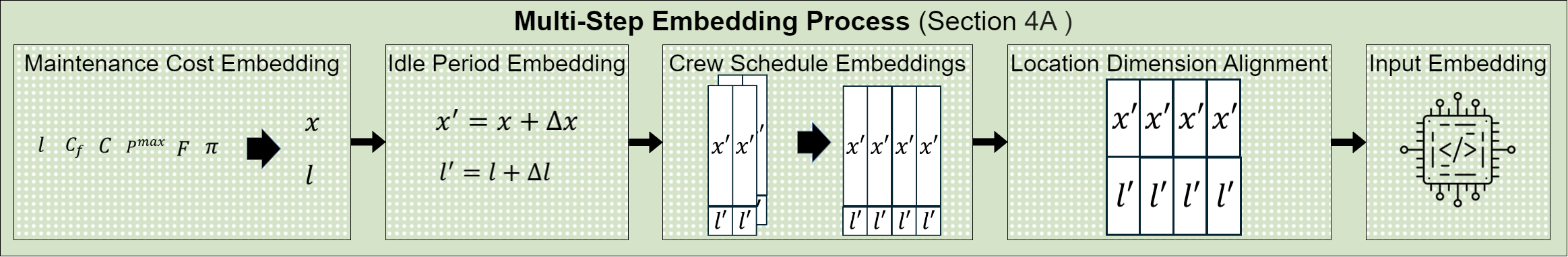}
    \caption{Multi Step Embedding in AttenCOpt.}
    \label{fig:step1}
\end{figure*}

In line with the MIP formulation, the proposed AttenCOpt also models the O\&M decisions for a set of turbines \( \mathcal{I} = \{1, \ldots, I\} \) over a planning horizon of \( T \) time steps with scenario set \( S \). Each turbine has scenario-specific temporal features for $s\in S$ such as the maximum production limit \( P^s_{i,t} \), and the electricity price \( \pi^s_{t} \) as well as the scenario agnostic predictive maintenance cost \( C_{i,t} \). Additionally, non-temporal features such as the failure time under scenario \( s \), \( F^s_{i} \), the unexpected failure cost \( C^f \), and the turbine location \( l_i \) are also considered. Together, these features form the input set for the problem as illustrated in Figure~\ref{fig:step1}. 

The following two sections define the training model which helps embed the MIP model into the MHA framework to ensure that the trained model effectively captures the spatial distribution of turbine locations, scenario-specific variations in production limits, electricity prices, and predictive maintenance costs. We discuss the embedding mechanism for integrating the objective function of Problem ~\ref{obj} into an MHA framework in Section \ref{subsec: mha_embed}. Further, the formulation of the MHA based reward functions are described in Section \ref{subsec:mha_reward}. 

\subsection{Multi-Step Embeddings for Wind Farm Maintenance} \label{subsec: mha_embed}
In order to successfully solve Problem \ref{obj} using an MHA framework, we need to obtain MHA compatible embeddings for each of the five maintenance input features. The methodologies to generate these embeddings are discussed below.

\subsubsection{Maintenance Cost Embeddings}
We consider preventive maintenance costs, failure cost, and failure scenarios as foundational input features for maintenance scheduling that are leveraged to yield MHA compatible embeddings. Specifically, for turbine \( i \) in period \( t \), we consider the maintenance cost as well as the profit denoted by $\quad \frac{1}{|S|} [ \sum_{s=1}^{S} \sum_{i=1}^{I} \sum_{t=1}^{F^s_i-1} \left( C_{i,t} + \pi^s_{t} \cdot P^s_{i,t} \right) \cdot m_{i,t} + 
\sum_{s=1}^{S} \sum_{i=1}^{I} \sum_{t=F^s_i}^{T} \left( C^f + \sum_{l=F^s_i}^t \pi^s_{l} \cdot P^s_{i,l} \right) \cdot m_{i,t}]$. When maintenance is performed before a failure event, the total cost includes the sum of the dynamic maintenance cost and the production losses due to downtime during the maintenance period. Conversely, when a failure occurs before maintenance, the total cost includes the sum of the unexpected failure cost and the production losses incurred from the time of failure until the maintenance is completed. As a result the total maintenance scheduling cost \( x_{i,t} \) for $i \in I$, $t \in T$ is denoted by Equation \ref{eq:Input Features Calculation}.
\begin{equation}
\begin{aligned}
\label{eq:Input Features Calculation}
x_{i,t} &= \frac{1}{|S|} \sum_{s \in S} \Bigg[ \mathbb{I}(t < F^s_{i}) \left( C_{i,t} + \pi^s_{t} P^s_{i,t} \right) \\
&\quad + \mathbb{I}(t \geq F^s_{i}) ( C^f + \sum_{l=F^s_{i}}^{t} \pi^s_{l} P^s_{i,l} ) \Bigg]
\end{aligned}
\end{equation}
The maintenance cost embedding can be used to inform the maintenance scheduling cost for the AttenCOpt for a fleet of turbines depending a given set of scenarios.

\subsubsection{Idle Period Embeddings}
A sequence-to-sequence driven MHA model necessitates scheduling exactly $M$ maintenance decisions in each period. However, the maintenance scheduling constraints encoded in Constraint ~\ref{eq:Maintenance Scheduling Constraints} in Problem \ref{obj} necessitates at most $M$ maintenance decisions per period. In order to resolve this inconsistency, we make the assumption that maintenance crews are located at the depot or base station when no maintenance activity is scheduled. As a result, we consider a set of idle turbines (\( I' \)) such that Equations \ref{eq:idle1} and \ref{eq:idle2} hold.
\begin{gather}
x_{i,t} = 0, \quad \forall t \in I', \; \forall i \label{eq:idle1}\\
l_i = 0, \quad \forall t \in I', \; \forall i \label{eq:idle2}
\end{gather}
Equations \ref{eq:idle1} and \ref{eq:idle2} ensure that each turbine in set $I'$ incurs no maintenance cost while being located at the depot or base station respectively. The implication of Equations \ref{eq:idle1} and \ref{eq:idle2} is that crew locations coincide with the depot whenever a maintenance of the idle turbines is scheduled. As a result, the set of idle turbines can effectively help MHA models relax the necessary constraint (of $M$ decisions per period) and  seamlessly help integrate Constraint ~\ref{eq:Maintenance Scheduling Constraints} of Problem \ref{obj} into the AttenCOpt. 

\subsubsection{Crew Schedule Embeddings}
\label{sec:expand}
The third step, \textit{Scheduling Multiple Maintenance Teams}, addresses the model's capability to perform multiple maintenance operations within each period, denoted by \( M \) (constraint~\ref{eq:Maintenance Scheduling Constraints}). To effectively represent this capacity, we duplicate the input features \( x_{t,i} \), corresponding to the maintenance scheduling cost for asset \( i \) at time period \( t \), based on the number of allowable maintenance operations per period.

\begin{figure}[h]
     \centering
     \includegraphics[width=0.6\linewidth]{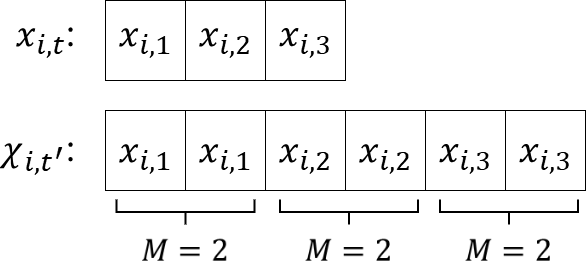} 
     \caption{Scheduling Multiple Maintenance Teams}
     \label{fig:expansion}
 \end{figure}

Constraint ~\ref{eq:Maintenance Scheduling Constraints} in Problem \ref{obj} opens up possibilities for multiple maintenance operations in each period. Therefore, we create embedding $\chi_i$ to represent $x_i$ for each turbine $i$ such that $\chi^m_{i,t'} = x_{i,t}$, $\forall m \in M$. As a result, the embedding $\chi_i$ will be of dimension $T \times M$. Figure \ref{fig:expansion} represents an illustrative example for a turbine $i$ with $M=2$.

\subsubsection{Location Dimension Alignment}

The fourth step, \textit{Location Dimension Alignment}, ensures the consistency of location data dimensions with the maintenance decision cost features. Initially, the location of each turbine is in \(\mathbb{R}(I + I')\), while \(\chi\) is in \(\mathbb{R}(I + I') \times T \times M\). To align these dimensions, the location vector \(l'\) is repeated \(T \times M\) times, matching the dimensionality of \(\chi\). This repetition integrates the location data with other input features, ensuring consistency across the dataset.

\subsubsection{Input Embedding}

After Step 4, the "Input Embedding" layer transforms inputs \((\chi, l')\) into a unified \( D \)-dimensional space, capturing patterns and dependencies for efficient maintenance scheduling and enhancing model performance.

\subsection{Reward Function Formulation}\label{subsec:mha_reward}

The cost of maintaining turbine \( v \) at time \( t \), after previously maintaining turbine \( u \), is represented by the cost function \( f_c : \chi_{u,t-1} \times \chi_{v,t} \rightarrow \mathbb{R} \).
Our objective is to determine a sequence of turbines \( Y \) of length \( T_m \in (0, T \times M] \) that satisfies the constraint~\ref{eq:Maintenance Constraints} of the wind farm maintenance scheduling problem (\( C(Y) \)), where \( T \times M\) accounts for the extended time periods due to the repetition of maintenance times per period (section~\ref{sec:expand}). We aim to minimize:
\begin{equation}
\label{eq:Objective Function}
P_{\text{obj}}[Y|I] = \sum_{t=1}^{T \times M} f_{c}(\chi(t-1, u), \chi(t, v))
\end{equation}

where \( y,v \in I \) are the turbines selected from \( I \) at time step \( t \). The function \( C \) checks whether the sequence \( Y \) satisfies constraint~\ref{eq:Maintenance Constraints}.

The function \(f_{c}(\chi(t-1, u), \chi(t, v)) \) is defined as:
\begin{equation}
\label{eq:cost_function}
\begin{aligned}
f_{c}(\chi(t-1, u), \chi(t, v)) &= \chi(t, v) + \mathbb{I}(\delta_t = 1) \cdot \Delta
\end{aligned}
\end{equation}
where \( \chi(t, v) \) is the maintenance decision cost of turbine \( v \) at time \( t \), \( \delta_{t} \) is 1 if turbine \( y \) is located at different location than turbine \(u\), and \(\Delta\) is the location change cost. This function captures the additional cost due to location changes between consecutive maintenance operations. This reward function is the same as the objective function defined in equation~\ref{obj}

\section{Attention model}
The objective of the wind turbine maintenance problem is to find a policy \( \pi \), that generates a sequence \( Y \) while minimizing \( P_{\text{obj}} \) and satisfying \( C(Y) \) for a given problem instance \( s \). The MHA-driven  maintenance scheduling (AttenCOpt) framework is presented in Figure \ref{fig: Main_fig}. The AttenCOpt essentially comprises an encoder and decoder model that defines a stochastic policy \( p(\pi | s) \) and is parametrized by $\theta$ as represented in Equation \ref{eq:pi cost} 

\begin{equation}
\label{eq:pi cost}
p_{\theta}(\pi | s) = \prod_{t=1}^{n} p_{\theta}(\pi_t | s, \pi_{1:t-1}).
\end{equation}

\subsection{Encoder}

To solve the wind turbine maintenance scheduling problem (Equation \ref{eq:pi cost}), we develop leverage the Graph Temporal Attention (GTA) Neural Network architecture from \cite{gunarathna2022solving}. The encoder is based on the transformer architecture introduced by \cite{vaswani2017attention}, as shown in Figure \ref{fig: Main_fig}. The temporal encoder receives input from a multi-step embedding and consists of parallel spatial and temporal attention layers. The outputs (\( H_S^{(1)}, H_T^{(1)} \)) from these layers are combined, and the first encoding layer is illustrated in Figure \ref{fig: Main_fig}. Additional layers can be stacked, using the outputs \( H^{(l)} \) given hidden state of turbine \( i \) at time \( t \) in layer \( l \) denoted by \( h_{i,t}^{(l)} \). The encoder in turn comprises of three components.

\subsubsection{Spatial Attention Layer} 
The spatial attention layer encodes dependencies between turbines at a given time step. Utilizing the MHA mechanism ~\cite{vaswani2017attention}, it processes a sequence of turbine features to capture inter-turbine relationships. The output of this layer is denoted as \( H_S^{(l+1)} \). The detailed mathematical operations for computing the spatial attention are provided in Appendix~B.

\subsubsection{Temporal Attention Layer}

Due to the dynamic nature of maintenance decision problems, dependencies exist between turbine features across different time steps in wind farm turbine maintenance problems. The temporal attention layer encodes these temporal dependencies by capturing relevant information from various time steps for each turbine. The output of this layer is denoted by \( H_T^{(l+1)} \). The mathematical formulation is presented in Appendix~C.

\subsubsection{Multi-Head Attention Framework}

We apply a linear transformation on the concatenation of the spatial  and temporal representations (denoted by (\( H_{T}^{(l+1)} \)), (\( H_S^{(l+1)} \)) respectively) which serves as an input for the sigmoid activation function.

\begin{equation}
\label{eq:MHAT}
H^{l+1} = \sigma(wI \cdot (H_S^{l+1} ||  H_T^{l+1})) 
\end{equation}

The integration layer output is designed to account for the influence of surrounding turbines as well as the temporal variations—both past and future—of each turbine.

\subsection{Decoder}
The decoder operates with a feedback loop that starts from the spatio-temporal output of the encoder's last layer \( H^{(L)} \) and the partially decoded solution up to time \( t \). The decoder architecture is illustrated in Figure \ref{fig: Main_fig}. At each time step, the temporal pointer generates a turbine embedding—excluding turbines that have already been maintained—and creates a context embedding based on the current solution. These embeddings are then combined using a multi-head attention layer. A log-probability layer determines iteratively helps select the next turbine.

\subsubsection{Temporal Pointer}

To capture dynamic dependencies at each time step, the temporal pointer computes attention by focusing on the most relevant parts of the embedded output from the encoder. Specifically, it slices the encoder embedding at each decoding time step and computes attention weights using multi-head attention. We denote this slice of the encoder embedding at time \( t \) as \( H_{D,t} \). This mechanism allows turbine representations to be updated at every time step based on newly computed values. The mathematical operations for the temporal pointer are detailed in Appendix~E.

\subsubsection{Context Embedding}
To define the context of the problem at each decoding time step, we introduce a problem-dependent context embedding \( H_{C} \) tailored to the wind farm maintenance scheduling problem. This embedding encompasses the last selected turbine and the aggregated problem embedding, computed by summing all turbine features at time \( t \). The detailed mathematical formulation of \( H_{C} \) is presented in Appendix~F.

\subsubsection{Masked Multi-head and Log Probability Layer}
After obtaining \( H_{C,t} \), a masked MHA layer integrates the context-specific embedding with the current problem representation. The query, key, and value weights are applied as in the same manner as in the encoder. The resulting embeddings are then used to compute the log probabilities for each turbine. The detailed mathematical formulation is shown in Appendix G.

\subsubsection{Reformulating linear constraint via masking procedure}
To ensure that the solutions comply with constraints, we employ a masking method in the decoder that aligns linear optimization constraints with their projections through conditional masking. This method ensures compliance with constraints such as avoiding simultaneous maintenance on certain turbines and prevents the selection of turbines that have already been maintained. All constraints or sets of constraints can be represented by \( C \). For example, in this problem, \( C \) denotes the constraint specified in Equation~\ref{eq:Maintenance Constraints}. The algorithm for the masking procedure is detailed in Algorithm~\ref{alg:masking_procedure}. In summary, if selecting a turbine for maintenance in the next period violates a constraint, a large negative value can be artificially added to its probability, effectively preventing it from being chosen as the next turbine to maintain.

\begin{algorithm}
\caption{Masking Procedure}
\label{alg:masking_procedure}
\begin{algorithmic}[1]
\STATE \textbf{Input:} Sequence of turbines \( Y \), Constraint matrix \( C \)
\STATE \textbf{Output:} Masked probabilities \( P_t \)
\FOR{each time step \( t \)}
    \STATE Initialize mask \( M \) as a zero vector 
    \FOR{each turbine \( y_t \) in \( Y \)}
        \STATE Set \( M[y_t] = -\infty \) \COMMENT{Mask already maintained turbines}
    \ENDFOR
    \FOR{each constraint \( c \) in \( C \)}
        \IF{constraint \( c \) is violated at time \( t \)}
            \STATE Update mask \( M \) according to \( c \)
        \ENDIF
    \ENDFOR
    \STATE Compute masked probabilities \( P_t = \text{softmax}(\gamma_t + M) \)
\ENDFOR
\end{algorithmic}
\end{algorithm}

\subsection{Reinforcement Learning driven AttenCOpt}
Our model is trained using a reinforcement learning approach, where the probability distribution \( p_{\theta}(\pi | s) \) defines a solution \( \pi | s \) for a given instance \( s \). The loss function \( \mathcal{L}(\theta | s) \) is defined as the expected cost, which combines dynamic maintenance costs and relocation of the maintenance team. The training process uses the gradient descent method, employing a reinforcement learning gradient estimator as described in \cite{williams1992simple}, with a baseline function \( b(s) \):

\begin{equation}
\label{eq:mha_decoder}
\nabla_{\theta} \mathcal{L}(\theta | s) = \mathbb{E}_{p_{\theta}(\pi | s)} \left[ (L(\pi) - b(s)) \nabla_{\theta} \log p_{\theta}(\pi | s) \right]
\end{equation}

We use a baseline \( b(s) \) to reduce gradient variance and accelerate learning by providing a reference point for the cost function, stabilizing the updates during training. A rollout baseline is chosen for this purpose; it is calculated by simulating the policy’s decisions multiple times (rollouts) to obtain an average cost estimate for the solution. For each state \( s \), the baseline function \( b(s) \) is the average cost computed from these rollouts, which approximates the expected cost of following the policy from state \( s \). The complete training algorithm is presented in Algorithm~\ref{alg:training}.

\begin{algorithm}[ht]
\caption{Training AttenCOpt with Reinforcement Learning}
\label{alg:training}
\begin{algorithmic}[1]
\STATE Initialize parameters $\theta$ for the AttenCOpt
\STATE Define the loss function \( \mathcal{L}(\theta | s) \)
\FOR{each iteration}
    \STATE Sample a set of instances \( s \) from the training set
    \FOR{each instance \( s \)}
        \STATE Sample a solution \( \pi \) from the policy \( p_{\theta}(\pi | s) \)
        \STATE Compute the cost \( L(\pi) \)
    \ENDFOR
    \IF{using baseline function}
        \STATE Compute the baseline function \( b(s) \)
        \STATE Compute the gradient using the reinforcement learning gradient estimator:
        \[
        \nabla_{\theta} \mathcal{L}(\theta | s) = \mathbb{E}_{p_{\theta}(\pi | s)} \left[ (L(\pi) - b(s)) \nabla_{\theta} \log p_{\theta}(\pi | s) \right]
        \]
    \ELSE
        \STATE Compute the gradient without baseline:
        \[
        \nabla_{\theta} \mathcal{L}(\theta | s) = \mathbb{E}_{p_{\theta}(\pi | s)} \left[ L(\pi) \nabla_{\theta} \log p_{\theta}(\pi | s) \right]
        \]
    \ENDIF
    \STATE Update the parameters $\theta$ using gradient descent:
    \[
    \theta \leftarrow \theta - \alpha \nabla_{\theta} \mathcal{L}(\theta | s)
    \]
\ENDFOR
\end{algorithmic}
\end{algorithm}

\section{Experiments}

We perform an extensive set of experiments to showcase the performance of the proposed AttenCOpt in different problem settings. 
The first part of our analysis is a comparison of the performance of AttenCOpt with the results obtained by using MIP model as a benchmark solved by Gurobi. For these set of experiments we will compare the methods using solution time and quality metrics. The second set of experiments will focus on how the proposed AttenCOpt trained on a specific case can be transferred across other cases. To this end, we will train and test the proposed model considering different cases, and showcase the transfer learning capability across different cases. 

We now specify the hyperparameters utilized in our experimental setup. In the case of Multi-head attention framework, we set \( N = 3 \) layers for the encoder, therefore, we select \( D^h = 128 \) as a hidden dimension. This configuration provides a good balance between quality and computational complexity. Using this scheme, the learning proceeds for 100 epochs. Each time, 25,600 instances constitute the source of problems, with a batch size of 32 instances. A learning rate of \( 1 \times 10^{-4} \) is applied, as higher learning rates showed instability during learning. The Adam optimizer is utilized for training. Our experiments were run on the supercomputer equipped with an Intel processor and 2 NVIDIA Quadro RTX6000 GPUs. 

\begin{figure}[htbp]
    \centering
    \begin{minipage}{.9\columnwidth}
        \centering
        \includegraphics[width=1\linewidth]{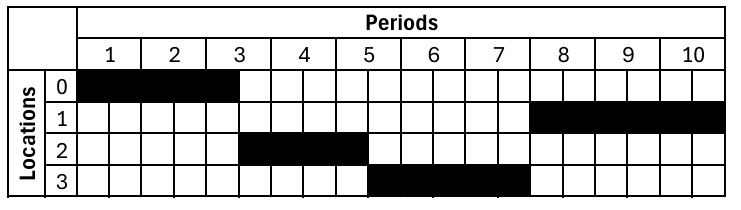} 
        \caption*{(a)}
    \end{minipage}%
    \vspace{1em}
    \begin{minipage}{.9\columnwidth}
        \centering
        \includegraphics[width=1\linewidth]{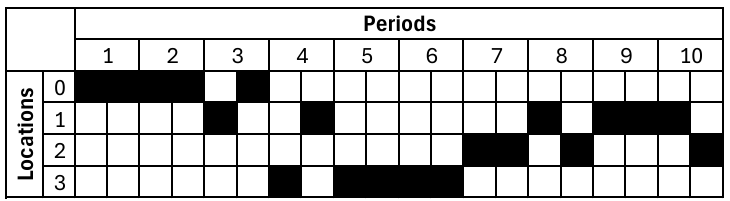}
        \caption*{(b)}
    \end{minipage}
    \caption{Demonstration of AttenCOpt results for maintenance scheduling of 20 turbines over 10 periods across 4 locations.}
    \label{fig:maintenance_scheduling}
\end{figure}

Figure \ref{fig:maintenance_scheduling} showcase a sample maintenance scheduling result obtained for 20 turbines over 10 periods across 4 locations, with '0' defined as the depot. Each period allows up to two maintenance activities.

\begin{itemize}
    \item \emph{Figure 4(a)} shows scheduling with positive location change costs (in objective function \(\Delta \neq 0\)), where maintenance activities are grouped to minimize costs.
    \item \emph{Figure 4(b)} shows scheduling without location change costs (in objective function \(\Delta = 0\)), resulting in dispersed activities across locations and periods, potentially increasing complexity and inefficiencies.
\end{itemize}

These contrasting scenarios demonstrate the impact of location change costs on optimizing maintenance operations, and also showcases how the proposed AttenCOpt adapts to visit costs to perform opportunistic maintenance. Following this sample demonstration, we will shift our focus to large-scale experiments.

\subsection{Comparative performance of the proposed and benchmark methods as a function of solution time and quality}

This section compares the performance of the AttenCOpt with MIP model solving with Gurobi, across various maintenance cases (see Table \ref{tab:descriptive_statistics}).
Each case differs by the number of turbines, locations, maintenance operations per period, and total periods. Gurobi solves small-size cases (up to 30 turbines and 20 stochastic scenarios) within 1 hour, but for larger cases, it becomes intractable, with computation times exceeding 7500 seconds.

\begin{table}[ht]
\centering
\footnotesize 
\setlength{\tabcolsep}{5pt} 
\renewcommand{\arraystretch}{1.1} 
\begin{tabular}{ccccccc}
\toprule
 & \textbf{Case 1} & \textbf{Case 2} & \textbf{Case 3} & \textbf{Case 4} & \textbf{Case 5}  \\
\midrule
\textbf{\# of Turbines} & 15 & 25 & 30 & 40 & 50  \\ 
\textbf{\# of Locations} & 4 & 4 & 4 & 4 & 4  \\ 
\textbf{Maint. per period} & 2 & 2 & 2 & 2 & 2  \\ 
\textbf{\# of Periods} & 10 & 15 & 20 & 25 & 30  \\ 
\midrule
\textbf{Mean} & 2.67 & 103.70 & 953 & 3268.18 & 3600  \\ 
\textbf{Median} & 2 & 70 & 510 & 3600 & 3600  \\ 
\textbf{Q1} & 1 & 60 & 340 & 3600 & 3600  \\ 
\textbf{Q3} & 5.00 & 120 & 1080 & 3600 & 3600  \\ 
\textbf{\% Solved} & 100\% & 100\% & 100\% & 20\% & 0\%  \\
\bottomrule
\end{tabular}
\caption{Solution times for MIP model across different cases. Corresponding solution times for AttenCOpt was excluded, as it solves all cases under one second.}
\label{tab:descriptive_statistics}
\end{table}

Table \ref{tab:descriptive_statistics} shows solution times for five cases using MIP formulation, including the mean, median, first quartile (Q1), and third quartile (Q3). \% Solved metric reports the ratio of problems that were solved by Gurobi within the specified time limit of 3600 seconds. The mean and median solving times increase significantly from Case 1 to Case 5, with Cases 4 and 5 reaching the 3600-second time cap.

The percentage of problems solved within 3600 seconds shows Gurobi's performance dropping from 100\% in Cases 1 to 3 to 20\% in Case 4 and none in Case 5. The AttenCOpt's solution times are excluded from this table, as our method solves all cases in under one second, demonstrating superior performance.

Given the superiority of the proposed method in solution times, the next metric to consider would be the quality of the provided solution. Table \ref{tab:maintenance_cases} showcases the solution quality of the proposed model as a function of the optimality gap. Our results show that the AttenCOpt can find high-quality solutions in real-time. For the largest case with 50 turbines, 4 locations, and 30 periods, our model solves the problem in half a second, and the solution gap between our model and Gurobi remains minimal, with a maximum difference of 1.57\% in the largest case.

\begin{table}[ht]
\centering
\footnotesize 
\setlength{\tabcolsep}{5pt} 
\renewcommand{\arraystretch}{1.1} 
\begin{tabular}{ccccccc}
\toprule
 & \textbf{Case 1} & \textbf{Case 2} & \textbf{Case 3} & \textbf{Case 4} & \textbf{Case 5}  \\
\midrule
\textbf{Mean Gap} & 0.13\% & 0.26\% & 0.24\% & 0.73\% & 1.57\%  \\ 
\midrule
\textbf{Q1 of Gap} & 0.08\% & 0.17\% & 0.17\% & 0.28\% & 0.35\%  \\ 
\textbf{Median of Gap} & 0.12\% & 0.19\% & 0.20\% & 0.85\% & 1.39\%  \\ 
\textbf{Q3 of Gap} & 0.15\% & 0.22\% & 0.24\% & 1.04\% & 2.70\%  \\ 
\textbf{Std of Gap} & 0.02\% & 0.29\% & 0.13\% & 0.51\% & 1.56\%  \\ 
\bottomrule
\end{tabular}
\caption{Solution quality of proposed multi-head attention framework across different cases as a function of the optimality gap.}
\label{tab:maintenance_cases}
\end{table}

\subsection{Transfer learning across different problem cases}

\begin{figure}
    \centering
    \subfloat[Model Trained on Case 1 Data]{%
        \includegraphics[width=0.75\linewidth]{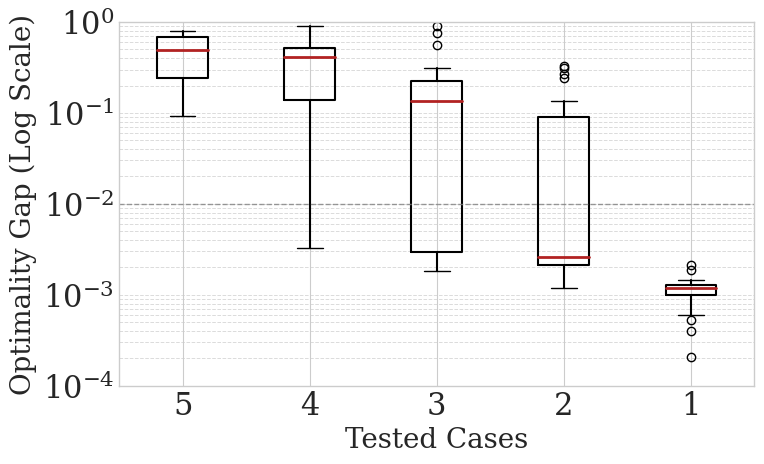}
        \label{fig:case_1}}
    \vspace{0.1cm} 
    \subfloat[Model Trained on Case 2 Data]{%
        \includegraphics[width=0.75\linewidth]{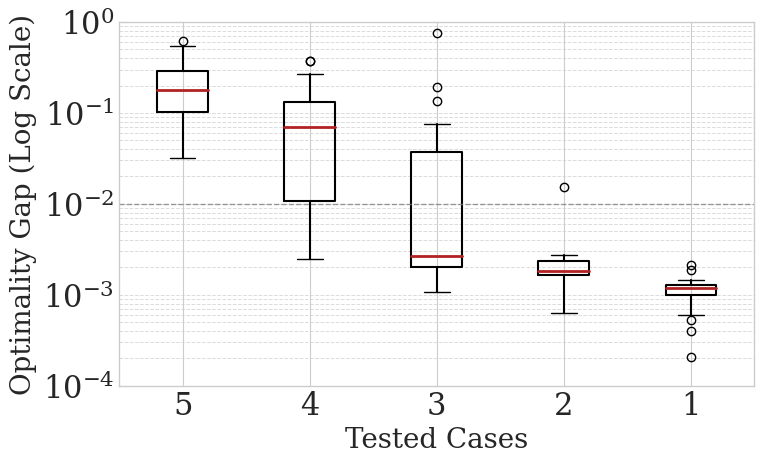}
        \label{fig:case_2}}
    \vspace{0.1cm}
    \subfloat[Model Trained on Case 3 Data]{%
        \includegraphics[width=0.75\linewidth]{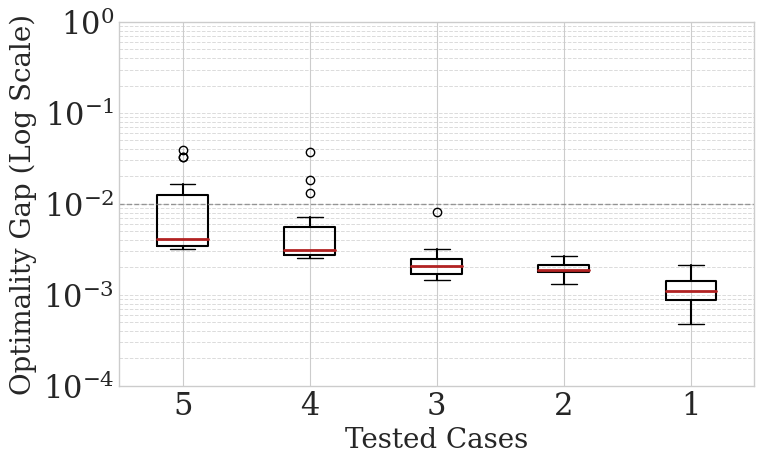}
        \label{fig:case_3}}
    \vspace{0.1cm}
    \subfloat[Model Trained on Case 4 Data]{%
        \includegraphics[width=0.75\linewidth]{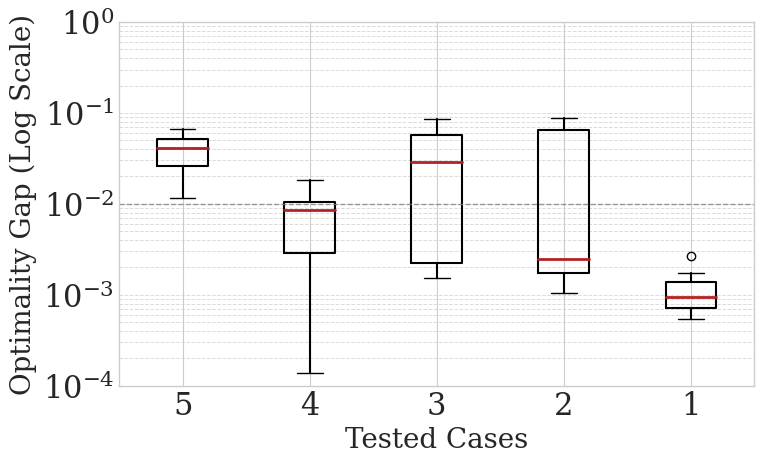}
        \label{fig:case_4}}
    \vspace{0.1cm}
    \subfloat[Model Trained on Case 5 Data]{%
        \includegraphics[width=0.75\linewidth]{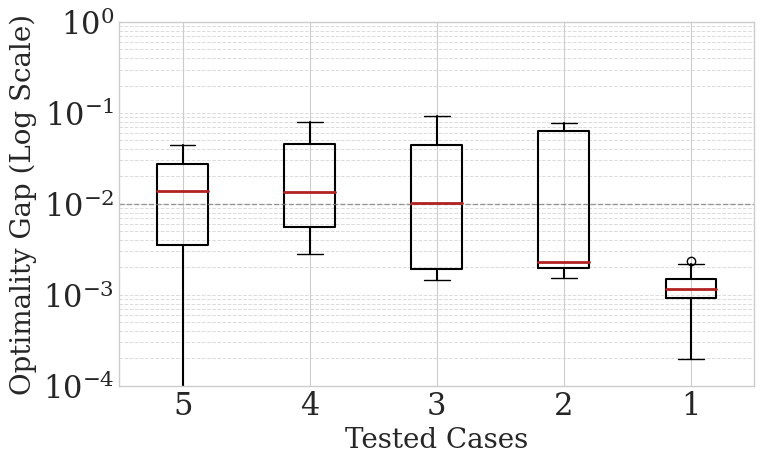}
        \label{fig:case_5}}
    
    \caption{Performance of different cases when trained on respective cases. Each subfigure shows the performance of various cases when trained on the corresponding case indicated.}
    \label{fig:cases}
\end{figure}

We evaluate the generalization capabilities of the AttenCOpt by evaluating the transfer learning capability. To this end, we train our MHA model  on one case, and showcase the quality of the solutions when it is tested on different cases. Figure~\ref{fig:cases} presents the model's performance for models trained and tested across cases 1 through 5. Each subfigure in the Figure~\ref{fig:cases} contains a average gap (in log scale) and the distribution of gap values using a box chart.

Figure~\ref{fig:cases} shows some key insights on our MHA model:
\begin{itemize}
\item \textit{MHA showcases significant transfer learning capability:} MHA models trained and tested in different cases typically result in optimality gaps less than 10\%, showcasing that there is significant transfer learning across cases.

\item \emph{Transfer learning capability is highest for models trained in large-cases:} In our experiments, we show that the models trained in large cases, e.g. case 5, typically result in very high solution quality for other tested cases. For instance, Case 5 yields a mean gap of approximately 0.123\% to 2.81\% across all the cases. This is likely a result of smaller cases being a subset of the large cases, dynamics of which are inherently learned during the training of large cases.
\end{itemize} 

These set of observations have interesting managerial implications. Firstly, we observe that the proposed AttenCOpt do not need to be trained for every problem case to be considered. In fact, our second observation suggests that the best approach is the following: (i) identify the number of turbines and time periods in the largest case that needs to be solved, and (ii) train the AttenCOpt for this large case. 
Subsequently, they can apply this model to solve problems of varying sizes with high efficiency and stability. This approach not only ensures robustness across different scales but also significantly reduces the need for retraining models for each specific case, thereby optimizing both time and computational resources.

\section{Conclusion}
In this work, we developed an MHA-based decision making framework for solving large-scale wind farm O\&M problems which we refer to as \textit{AttenCOpt}. The proposed model (i) significantly reduces the solution time from hours to seconds, (ii) guarantees feasibility of the proposed solutions considering complex constraints that are omnipresent in wind farm O\&M, and (iii) results in significant solution quality compared to the conventional MIP formulations. 
The development of the proposed MHA framework incorporates a number of modeling innovations that allows explicit embedding of objective function, constrains and problem definitions across different components of the MHA model. The objective function and part of the constraints were reformulated within the reward function, while the remainder of the constraints were incorporated using the masking procedure. The proposed transition from MIP to MHA model is not problem specific, and can be implemented for a range of O\&M scheduling problems in wind farms with small alterations.

The proposed Multi-head Attention Framework advances wind farm maintenance scheduling by integrating ML techniques and attention mechanisms with MIP models. It offers a scalable, efficient, and adaptive solution, addressing current challenges and enabling future innovations in predictive maintenance and operational optimization in renewable energy. Future research may extend its application to other domains and further improve its performance.

\bibliographystyle{unsrt}
\bibliography{refe}

\newpage
\pagenumbering{arabic}
\setcounter{page}{1} 
\appendix
\label{appendix:proof1}

\subsection*{Appendix A: Proof of Equivalence of Objective Functions}

Consider the original objective function:

\begin{equation}
\begin{aligned}
\text{Maximize} \quad \frac{1}{|S|} \biggl[ &\sum_{s=1}^{S} \sum_{i=1}^{I} \sum_{t=1}^{T}  \left(\pi^s_{t} \cdot y^s_{i,t} \right) - \\
&\biggl[ \sum_{s=1}^{S} \sum_{i=1}^{I} \sum_{t=1}^{F^s_i-1}  \left(C_{i,t} \cdot m_{i,t}\right) + \sum_{s=1}^{S} \sum_{i=1}^{I} \sum_{t=F^s_i}^{T} \\
&\left( C^f \cdot m_{i,t} \right) \biggr] \biggr] - \sum_{t=1}^T \delta_t \cdot \Delta \\
&\text{S.t:}~(
\ref{eq:production_constraint}-\ref{eq:domain})
\end{aligned}
\end{equation}

This can be rewritten in terms of a subproblem \( q(m) \):

\begin{equation}
\begin{aligned}
\text{Maximize} \quad q(m) + \frac{1}{|S|} \biggl[ &- \sum_{s=1}^{S} \sum_{i=1}^{I} \sum_{t=1}^{F^s_i-1} C_{it} m_{i,t} \\
&- \sum_{s=1}^{S} \sum_{i=1}^{I} \sum_{t=F^s_i}^{T} C^f m_{i,t} \biggr] 
- \sum_{t=1}^T \delta_t \Delta\\
&\text{S.t:}~(
\ref{eq:Maintenance Scheduling Constraints}-\ref{eq:domain})
\end{aligned}
\end{equation}

where

\begin{align}
&q(m) = \max \frac{1}{|S|} \sum_{s=1}^{S} \sum_{i=1}^{I} \sum_{t=1}^{T} \pi^s_{t} y^s_{i,t} \\
&\text{S.t: } (\ref{eq:production_constraint} - \ref{eq:ailalibility_constraint2})
\end{align}

With constraints (\ref{eq:production_constraint} - \ref{eq:ailalibility_constraint2}), there is only one constraint for every time period \(t\), so for any fixed \(m_{i,t}\) the constraint (\ref{eq:production_constraint} - \ref{eq:ailalibility_constraint2}) can be be written as 

\begin{equation}
y^s_{i,t} \leq G^s_{i,t} \quad \text{for all } t \in T
\end{equation}

and

\begin{equation}
G^s_{i,t} = 
\begin{cases}
P^s_{i,t} (1 - m_{i,t}) & \text{for } t \leq F^s_i-1 \\
P^s_{i,t} \sum_{l=1}^{t-1} m_{i,l} & \text{for } t \geq F^s_i
\end{cases}
\end{equation}

Thus, the subproblem \( q(m) \) is maximized by choosing \( y^s_{i,t} = G^s_{i,t} \), giving \(\pi^s_{t} \geq 0\):

\begin{equation}
y^{s*}_{i,t} = G^{s}_{i,t}
\end{equation}

Substitute \( y^s_{i,t} = G^s_{i,t} \) into \( q(m) \):

\begin{equation}
q(m) = \frac{1}{|S|} \sum_{s=1}^{S} \sum_{i=1}^{I} \sum_{t=1}^{T} \pi^s_{t} G^s_{i,t}
\end{equation}

Now, substitute this expression back into the original objective function, and use the equality $\sum_{l=1}^{t-1} m_{i,l} = 1 - \sum_{l=t}^{T} m_{i,l}$ (which comes from Constraint~\ref{eq:Maintenance Constraints}):

\begin{equation}
\begin{aligned}
\text{Maximize} \quad & \frac{1}{|S|} \Biggl[ 
  \sum_{s=1}^{S} \sum_{i=1}^{I} \sum_{t=1}^{F_i^s - 1} \pi_{t}^s P_{i,t}^s (1 - m_{i,t}) \\
  & + \sum_{s=1}^{S} \sum_{i=1}^{I} \sum_{t=F_i^s}^{T}  \pi_{t}^s P_{i,t}^s (1-\sum_{l=t}^{T}  m_{i,l} )
\Biggr] \\
& - \Biggl[ 
  \sum_{s=1}^{S} \sum_{i=1}^{I} \sum_{t=1}^{F_i^s - 1} C_{i,t} m_{i,t} 
  + \sum_{s=1}^{S} \sum_{i=1}^{I} \sum_{t=F_i^s}^{T} C^f m_{i,t} 
\Biggr] \\
& - \sum_{t=1}^{T} \delta_t \Delta \\
= \; & \frac{1}{|S|} \Biggl[ 
  \sum_{s=1}^{S} \sum_{i=1}^{I} \sum_{t=1}^{T} \pi_{t}^s P_{i,t}^s 
  - \sum_{s=1}^{S} \sum_{i=1}^{I} \sum_{t=1}^{F_i^s - 1} m_{i,t} \pi_{t}^s P_{i,t}^s \\
  & - \sum_{s=1}^{S} \sum_{i=1}^{I} \sum_{t=F_i^s}^{T} \pi_{t}^s P_{i,t}^s  \sum_{l=t}^{T} m_{i,l} 
\Biggr] \\
& - \Biggl[ 
  \sum_{s=1}^{S} \sum_{i=1}^{I} \sum_{t=1}^{F_i^s - 1} C_{i,t} m_{i,t} 
  + \sum_{s=1}^{S} \sum_{i=1}^{I} \sum_{t=F_i^s}^{T} C^f m_{i,t} 
\Biggr] \\
& - \sum_{t=1}^{T} \delta_t \Delta
\end{aligned}
\end{equation}

Then, using the following equivalence:

\begin{equation}
 \sum_{s=1}^{S} \sum_{i=1}^{I} \sum_{t=F_i^s}^{T} \pi_{t}^s P_{i,t}^s \sum_{l=t}^{T} m_{i,l} =  \sum_{s=1}^{S} \sum_{i=1}^{I} 
\sum_{t=F_i^s}^{T} m_{i,t}
\sum_{l=F_i^s}^{t} \pi_{l}^s P_{i,l}^s    
\end{equation}

and eliminating the constant \( \sum_{s=1}^{S} \sum_{i=1}^{I} \sum_{t=1}^{T} (\pi^s_{t} \cdot P^s_{i,t}) \), we can obtain the following objective function:

\begin{equation}
\begin{aligned}
\text{Minimize} \quad \frac{1}{|S|} \biggl[ &\sum_{s=1}^{S} \sum_{i=1}^{I} \sum_{t=1}^{F^s_i-1} \left( C_{i,t} + \pi^s_{t} \cdot P^s_{i,t} \right) \cdot m_{i,t} + \\
&\sum_{s=1}^{S} \sum_{i=1}^{I} \sum_{t=F^s_i}^{T} \left( C^f + \sum_{l=F^s_i}^t \pi^s_{l} \cdot P^s_{i,l} \right) \cdot m_{i,t} \biggr] + \\
&\sum_{t=1}^T \delta_t \cdot \Delta.\\
&\text{S.t:}~(
\ref{eq:Maintenance Scheduling Constraints}-\ref{eq:domain})
\end{aligned}
\end{equation}

\subsection*{Appendix B: Mathematical Formulation of the Spatial Attention Layer}

The spatial attention layer employs the self-attention mechanism as introduced by \cite{vaswani2017attention}. For a given time step \( t \), the sequence of turbine features is represented as \( H_{S,t}^{(l)} = \{h_{1,t}^{(l)}, h_{2,t}^{(l)}, \ldots, h_{I,t}^{(l)}\} \in \mathbb{R}^{(I+I') \times D^h} \). The spatial attention uses standard parameters named query (\( w_{q}^{S(l)} \in \mathbb{R}^{D^h \times D^k} \)), key (\( w_{k}^{S(l)} \in \mathbb{R}^{D^h \times D^k} \)), and value (\( w_{v}^{S(l)} \in \mathbb{R}^{D^h \times D^v} \)).

First, we compute the following three projections:

\begin{equation}
q_{t}^{S(l)} = w_{q}^{S(l)} H_{S,t}^{(l)}, \quad k_{t}^{S(l)} = w_{k}^{S(l)} H_{S,t}^{(l)}, \quad v_{t}^{S(l)} = w_{v}^{S(l)} H_{S,t}^{(l)}.
\end{equation}

The scaled dot product between the query and key projections is computed, and a softmax layer is applied:

\begin{equation}
\alpha_t = \text{softmax}\left( \frac{q_{t}^{S(l)} (k_{t}^{S(l)})^T}{\sqrt{D^h}} \right).
\end{equation}

The output of the spatial attention for each turbine is:

\begin{equation}
H_{S,t}^{(l+1)} = \alpha_t v_{t}^{S(l)}.
\end{equation}

\subsection*{Appendix C: Mathematical Formulation of the Temporal Attention Layer}

The temporal attention layer encodes the dependencies between turbine features across different time steps. The input for this layer is the dynamic representation \( H^{(l)}_{T} \), and it computes the attention output \( H_{T,i}^{(l+1)} \) as follows:

\begin{equation}
q_{i}^{T(l)} = w_q^{T(l)} H_{T,i}^{(l)}, \quad k_{i}^{T(l)} = w_k^{T(l)} H_{T,i}^{(l)}, \quad v_{i}^{T(l)} = w_v^{T(l)} H_{T,i}^{(l)}.
\end{equation}

The scaled dot product between the query and key projections is computed:

\begin{equation}
\beta_i = \text{softmax}\left( \frac{q_i^{T(l)} (k_i^{T(l)})^T}{\sqrt{D^h}} \right).
\end{equation}

The output of the temporal attention layer is:

\begin{equation}
H_{T,i}^{(l+1)} = \beta_t v_{i}^{S(l)}.
\end{equation}

\subsection*{Appendix D: Multi-Head Attention framework}

The integration of spatial and temporal representations is achieved by concatenating \( H_S^{(l+1)} \) and transposed \( H_{T}^{(l+1)} \), followed by a linear transformation and sigmoid activation:

\begin{equation}
H^{(l+1)} = \sigma(wI \cdot (H_S^{(l+1)} \| H_{T}^{(l+1)})).
\end{equation}

The output \( H^{(l+1)} \) is then fed into the next layer in the temporal encoder.

\subsection*{Appendix E: Mathematical Formulation of the Temporal Pointer}

To dynamically compute attentions, we use the temporal pointer as defined in \cite{gunarathna2022solving}. The output \( H^{(L)} \) of the temporal encoder has dimensions \( \mathbb{R}^{(T \times M) \times (I+I') \times D^h} \). At each decoding time step, the encoder embedding is sliced, and multi-head attention weights are dynamically computed to focus on the most relevant parts of the embedded output. 

\subsection*{Appendix F: Mathematical Formulation of Context Embedding}

The context embedding \( H_{C} \) is used to identify the context at a decoding time step. For the wind farm maintenance scheduling problem, \( H_{C} \) contains the last selected turbine and the problem embedding computed by summing up all turbines at time \( t \):

\begin{equation}
H_{C,t} = \{ h^{(L)}_{y_t} \| H_{G,t} \}; \quad H_{G,t} = \sum_{i=0}^{I+I'} h^{(L)}_{i,t}
\end{equation}

\subsection*{Appendix G: Masked Multi-head Attention and Log Probability Layer}

The masked multi-head attention layer combines the context-specific embedding \( H_{C,t} \) with the current problem representation \( H_{D,t} \). Using the query, key, and value weights:

\begin{equation}
qC_t = w {C_q} H_{C,t}, \quad kC_t = w {C_k} H_{D,t}, \quad vC_t = w {C_v} H_{D,t}
\end{equation}

The attention weights and final embedding \( H_{D,t}^{(F)} \) are computed as:

\begin{equation}
H_{D,t}^{(F)} = \text{softmax}\left( \frac{qC_t^T kC_t}{\sqrt{D^h}} \right) \cdot vC_t
\end{equation}

Log probabilities for each node are computed using a weight vector \( wP \) and a tanh activation function:

\begin{equation}
\gamma_t = \tanh\left( H_{D,t}^{(F)} \cdot (wP \cdot H_{D,t})^T \right)
\end{equation}

A softmax layer computes the final probabilities for each turbine:

\begin{equation}
P_t = \text{softmax}(\gamma_t)
\end{equation}
\end{document}